\documentclass[12pt]{article}
\usepackage[english]{babel}
\usepackage{amssymb}
\usepackage[centertags]{amsmath}
\title{\textbf{On characteristic classes of  $Q$-manifolds}}
\author { S.L. Lyakhovich\footnote{Supported in part by  RFBR grant 05-01-00996.},
E.A. Mosman{}\footnote{Supported by RFBR grant 06-02-17352.}, and
A.A. Sharapov$^\dagger$}
\date{\footnotesize\textit{Department of Quantum Field Theory, Tomsk State University, Tomsk 634050, Russia }}
\begin{document}

\maketitle

\begin{abstract}
We define the notion of characteristic classes for supermanifolds
endowed with a ho\-mo\-lo\-gi\-cal vector field $Q$.  These take
values in the cohomology of the Lie derivative operator $L_Q$
acting on arbitrary tensor fields. We formulate a classification
theorem for intrinsic characteristic classes and give their
explicit description.
\end{abstract}

\noindent{\textbf{1.}} Let $M$ be a smooth supermanifold and
$\mathcal{T}(M)=\bigoplus_{n,m\in \mathbb{N}}
\mathcal{T}^{(n,m)}(M)$ be its tensor algebra; here
$\mathcal{T}^{(n,m)}(M)$ is the space of $n$-times contravariant
and  $m$-times covariant tensor fields on $M$. The elements of
$\mathcal{T}^{(1,1)}(M)$ are naturally identified with the
endomorphisms of the $C^{\infty}(M)$-module
$\mathcal{T}^{(1,0)}(M)$. The composition of two  endomorphisms
endows $\mathcal{T}^{(1,1)}(M)$ with the structure of an
associative algebra over $C^\infty(M)$. We will denote this
algebra $\mathcal{A}$. There is a natural trace on $\mathcal{A}$,
which is a $C^{\infty}(M)$-linear map $\mathrm{Str}:
\mathcal{A}\rightarrow C^\infty(M)$ vanishing on supercommutators.

\vspace{3mm}\noindent{\textbf{2.}} An odd vector field  $Q\in
\mathcal{T}^{(1,0)}(M)$ is said to be  \textit{homological}, if
\begin{equation}\label{QQ}
Q^2=\frac 12[Q,Q]=0\,.
\end{equation}
By definition \cite{Sch}, the pair $(M,Q)$ is called a
$Q$-\textit{manifold}.

The simplest example of a  $Q$-manifold is an odd tangent bundle
$\Pi TN$ (i.e., the tangent bundle of $N$ with reversed parity of
fibers). In this case, the supercommutative algebra of functions
$C^\infty(\Pi TN)$ is naturally isomorphic to the exterior algebra
of differential forms on $N$, with the de Rham  differential being
the (canonical) homological vector field on $\Pi TN$. A great
number of interesting examples of $Q$-manifolds is provided by Lie
algebroids \cite{V} and various gauge systems \cite{AKSZ, LS,
Sch}. For a recent discussion of homological vector fields in the
category of graded supermanifolds we refer the reader to \cite{M}.

Given a homological vector field $Q$, one can view
$\mathcal{T}(M)$ as a differential group  with the coboundary
operator
\begin{equation}\label{d}
\delta A=L_Q A\,,\qquad \forall A\in \mathcal{T}(M)\, .
\end{equation}
Here $L_Q$  is the Lie derivative w.r.t. $Q$. The property
$\delta^2=0$ follows from the identity $L^2_Q=L_{Q^2}=0$. Denote
by $H_Q(M)$ the group of $\delta$-cohomology. Since $\delta$
differentiates the tensor product, the group
$H_Q(M)=\bigoplus_{n,m\in \mathbb{N}} H^{(n,m)}_Q(M)$ inherits the
structure of a bigraded associative algebra over $\mathbb{R}$.

\vspace{3mm}\noindent{\textbf{3.}} Let $\nabla$ be a symmetric
connection on $M$ with the curvature tensor  $R_{XY}\in
\mathcal{A}$,
\begin{equation}
R_{XY}=[\nabla_X, \nabla_Y]-\nabla_{[ X, Y]} \,, \qquad \forall X,
Y \in \mathcal{T}^{(1,0)}(M)\,.
\end{equation}
Later on we will need the following relations characterizing the
geometry of $Q$-manifolds with symmetric connection:
\begin{equation}\label{rel}
\begin{array}{c}
\nabla_Q Q=0\,, \qquad \nabla_{Q}\Lambda = \frac{1}{2}R_{{Q} {Q}}
+ \Lambda^{2}\,,\qquad
 \nabla_{X}R_{Q Q} =2( R_{[
X,Q] Q} - \nabla_{{Q}}R_{Q X} )\,.
\end{array}
\end{equation}
Here  $X$ is an arbitrary vector field and $\Lambda \in
\mathcal{A}$ is an odd endomorphism  defined by the rule $\Lambda
(X)=\nabla_X Q$. With the tensor $\Lambda$, we have a simple
relation between the Lie and covariant derivatives of an
endomorphism $A\in \mathcal{A}$ w.r.t. the homological vector
field:
\begin{equation}\label{cov-L}
    \nabla_Q A=L_QA +[\Lambda, A]\,.
\end{equation}
A $Q$-manifold is said to be \textit{flat} if it admits a flat
connection.

\vspace{3mm}\noindent{\textbf{4.}} Under \textit{universal
cocycles} of a $Q$-manifold $M$ we understand  $\delta$-cocycles
$\mathcal{C}_\nabla [Q]\in \mathcal{T}(M)$ that are given, at each
coordinate chart,  by polynomials in the components of the
homological field $Q$, Cristoffel symbols of $\nabla$, and their
partial derivatives up to some finite order. The adjective
``universal'' emphasizes the fact that the closeness condition
$\delta (\mathcal{C}_\nabla[Q])=0$ is assumed to be satisfied by
virtue of equation (\ref{QQ}) without using a particular structure
of $Q$, $\nabla$, and $M$.

For example, the tensor powers of homological vector field
$Q^{\otimes n}\in \mathcal{T}^{(n,0)}(M)$ exhaust all  universal
cocycles that are independent of connection. A less trivial
example of universal cocycles is the function
\begin{equation}\label{P}
P_n = \mathrm{Str}((R_{QQ})^{2n})\in C^{\infty}(M)\,,\qquad n\in
\mathbb{N},
\end{equation}
which is nothing but the $n$th Pontrjagin's character of the
cotangent bundle $TM$ evaluated on the homological vector field
$Q$. (For the definition and discussion of the characteristic
classes of supervector bundles see \cite{VMP}.)

\vspace{3mm}\noindent\textbf{5.} The \textit{characteristic
classes of $Q$-manifolds} are, by definition, the elements of
$H_Q(M)$ that are represented by universal cocycles.

\vspace{3mm}\noindent\textbf{Theorem 1.}  \textit{The cohomology
classes of universal cocycles do not depend on the choice of
symmetric connection, and hence they are invariants of a
$Q$-manifold itself.}

\vspace{3mm}\noindent \textit{Proof}. Let
$\mathcal{C}_{\nabla_0}[Q]$ and $\mathcal{C}_{\nabla_1}[Q]$ be two
universal cocycles that differ only by the choice of connection.
Consider the direct product of $M$ and the linear superspace
$\mathbb{R}^{1,1}$ with one even coordinate $t$ and one odd
coordinate $\theta$. Equip the supermanifold
$\widetilde{M}=M\times \mathbb{R}^{1,1}$ with the homological
vector field  $ \widetilde{Q}=Q +\theta
\partial_t $ and  connection  $\widetilde{\nabla}=\nabla_t\oplus
\nabla'$, where  $\nabla_t=t\nabla_1+(1-t)\nabla_0$ is a
one-parameter  family of connections on $M$ and  $\nabla '$ is a
flat connection on $\mathbb{R}^{1,1}$. By universality,  the
tensor $\mathcal{C}_{\widetilde{\nabla}}[\widetilde{Q}]\in
\mathcal{T}(\widetilde{M})$ is closed w.r.t.
$\widetilde{\delta}=L_{\widetilde{Q}}$. Since $\theta^2=0$, one
can see that
$\mathcal{C}_{\widetilde{\nabla}}[\widetilde{Q}]=\mathcal{C}_{{\nabla_t}}[Q]+\theta
\Psi_{t}$, where $\Psi_{t}$ is some expression depending on  $Q$,
$\nabla_0$, $\nabla_1$, and $t$. The closeness condition
$\widetilde{\delta}(\mathcal{C}_{\widetilde{\nabla}}[\widetilde{Q}])=0$
is then equivalent to the following relations:
\begin{equation}\label{}
    \delta (\mathcal{C}_{{\nabla_t}}[Q])=0\,,\qquad \partial_t
    \mathcal{C}_{\nabla_t}[Q]=\delta \Psi_{t}\,.
\end{equation}
Integrating the second relation over $t$ from $0$ to $1$, we get
\begin{equation}\label{}
\mathcal{C}_{\nabla_1}[Q] - \mathcal{C}_{\nabla_0}[Q]=\delta
\int_0^1 dt \Psi_{t}\,.
\end{equation}
\nopagebreak[3] Thus, $\mathcal{C}_{\nabla_0}[Q]$ is cohomologous
to $\mathcal{C}_{\nabla_1}[Q]$. $\Box$

The first nontrivial series of charclasses of $Q$-manifolds was
proposed in \cite{LS}, the so-called \textit{principal series}.
The universal cocycles of this series involve the first covariant
derivatives of the homological vector field and have the following
form:
\begin{equation}\label{A}
C^{\infty}(M)\ni A_n
=\mathrm{Str}(\Lambda^{2n+1})+(\text{curvature dependent
terms})\,,\qquad \forall n\in \mathbb{N}\,.
\end{equation}
It was also shown that the class $A_0$ has a direct relationship
to one-loop anomalies in the BV quantization method of gauge
theories. In a particular case of  homological vector fields
corresponding to  Lie algebroids \cite{V}, formula (\ref{A})
reproduces the characteristic classes of Lie algebroids introduced
by Fernandes \cite{F}.

\vspace{3mm}\noindent{\textbf{6.}} In this note, we present two
infinite series of  universal cocycles, which essentially involve
the second covariant derivatives of the homological  vector field.
Together with the universal cocycles (\ref{A}) these new cocycles
generate and exhaust, in essence, all interesting
charac\-ter\-is\-tic classes of $Q$-manifolds. In order to write
them down in an explicit form we identify  $\mathcal{T}^{(1,
{}n+1)}(M)$ with $\mathcal{T}^{(0,n)}(M)\otimes \mathcal{A}$ and
treat the elements of the latter space as $n$-forms on $TM$ with
values in $\mathcal{A}$.

\vspace{2mm}\noindent \textbf{Lemma}. \textit{For any vector field
$X$, set $\Omega_X \equiv\nabla_X \Lambda -R_{XQ}\in \mathcal{A}$.
The tensor $\Omega \in \mathcal{T}^{(0,1)}(M)\otimes \mathcal{A}$
is a universal cocycle, i.e., $\delta \Omega =0$.}

\vspace{2mm}\noindent\textit{Proof}.  It follows from Rels.
(\ref{rel}) that
\begin{equation}
\nabla_Q \Omega_X =\Omega_{[ Q,X]}+[\Lambda,\Omega_X]\,,\qquad
\forall X\in \mathcal{T}^{(1,0)}(M)\,.
\end{equation}
Using successively the definition (\ref{d}), Rel. (\ref{cov-L}),
and the last identity, we find
\begin{equation}\label{}
    (\delta \Omega)_X=(L_Q\Omega)_X =
    L_Q(\Omega_X)-\Omega_{[Q,X]}=\nabla_Q \Omega_X -[\Lambda,
    \Omega_X]-\Omega_{[Q,X]}=0\,.
\end{equation}$\Box$.

Since the differential $\delta$ is compatible with contraction of
tensor indices, we have immediately

\vspace{2mm}\noindent \textbf{Corollary  1}. For any $n\in
\mathbb{N}$, define $B_n \in \mathcal{T}^{(0,n)}(M)\otimes
\mathcal{A}$ as
\begin{equation}\label{B}
B_n (X_1,...,X_n)= \Omega_{X_1} \Omega_{X_2}\cdots \Omega_{X_n}\,,
\quad \mathrm{for} \quad n>0\,,
\end{equation}
and  $B_0=1\in \mathcal{A}$. Then $\delta B_n=0$.

\vspace{2mm}\noindent \textbf{Corollary  2}. The  $n$-forms
\begin{equation}\label{C}
    C_n(X_1,...,X_n) =\mathrm{Str}(\Omega_{X_1} \Omega_{X_2}\cdots \Omega_{X_n})
\end{equation}
are $\delta$-closed and invariant, up to sign,  under the cyclic
permutation of their arguments:
\begin{equation}\label{}
    C_n(X_n,X_1,...,X_{n-1})=(-1)^{\varepsilon_1\varepsilon_2}
    C_n(X_1,...,X_n)\,.
\end{equation}
Here $\varepsilon_1 =\epsilon (X_n)+1$, $\varepsilon_2 =
\sum_{k=1}^{n-1}(\varepsilon (X_k)+1)$, and $\varepsilon(X_k)$
denotes the parity of the vector field $X_k$.

\vspace{2mm} We will refer to the cohomology classes of universal
cocycles (\ref{A}), (\ref{B}), and (\ref{C}) as the characteristic
classes of $A$, $B$, and $C$ series, respectively.

\vspace{3mm}\noindent{\textbf{7.}} We say that a characteristic
class $[\mathcal{C}_\nabla[Q]]\in H_Q(M)$ is \textit{intrinsic},
if it does not vanish identically upon setting the curvature of
$\nabla$ to zero. In other words, the intrinsic charclasses
survive on flat $Q$-manifolds. For example, all the characteristic
classes from the $A$, $B$, and $C$ series are intrinsic, while the
$\delta$-cohomology classes of (\ref{P}) are not. Clearly, the
intrinsic charclasses  constitute a subalgebra
$H^{\mathrm{int}}_Q(M)$ in $H_Q(M)$. As the next theorem shows,
this sub\-algeb\-ra admits a fairly simple description.

 \vspace{2mm}\noindent \textbf{Theorem 2}. \textit{The algebra $H^{\mathrm{int}}_Q(M)$
 is freely generated by the characteristic classes of A, B, and C series,
  together with the cohomology class
$[Q]\in H^{(1,0)}_Q(M)$ of the homological vector  field itself.}

 \vspace{2mm}\noindent\textit{Remark}. It might  be well to point out that
certain of the intrinsic charclasses may vanish for a particular
$Q$-manifold, e.g. by dimensional reasons. The theorem above
states just the absence of \textit{universal} nontrivial relations
between the generators of $H^\mathrm{int}_Q(M)$.

\vspace{2mm}

Let $V$ denote the typical fiber of the tangent bundle $TM$. The
proof of the theorem  is based on construction of a classifying
$Q$-map from a given \textit{flat} $Q$-manifold to the
infinite-dimensional, linear $Q$-manifold associated to the Lie
superalgebra $L_0(V)$ of formal vector fields on $V$ vanishing at
the origin. This reduces the problem to computation of stable
cohomologies of $L_0(V)$ with tensor coefficients. In the special
case that $V$ is an ordinary (even) linear space, the last problem
was completely solved in \cite{Fu} and the method of that paper
applies to the super case as well. The last but not least step
involves extension of the ``flat'' universal cocycles to arbitrary
(not necessarily flat) $Q$-manifolds.

 The details of the proof will be given elsewhere, along with various applications and
inter\-pre\-tat\-ions of the intrinsic characteristic classes.

\vspace{3mm} We wish to thank D.A. Leites for his useful comments
on the first version of the manuscript.

\makeatother

\begin{thebibliography}{1}


\bibitem{AKSZ} Alexandrov M., Kontsevich M., Schwarz A., Zaboronsky O., \textit{The
Geometry of the Master Equation and Topological Quantum Field
Theory}, Int. J. Mod. Phys. \textbf{A12} (1997) 1405-1430.

\bibitem{F} Fernandes R.L., \textit{Lie Algebroids, Holonomy and
Characteristic Classes}, Adv. in Math. \textbf{170}, N1 (2002),
119-179.

\bibitem{Fu} Fuks D.B., \textit{Stable cohomologies of a Lie algebra of formal vector fields with tensor
coefficients}, Funct. Anal. Appl.,  17 (1983), 295-301.



\bibitem{LS}Lyakhovich S.L., Sharapov A.A., \textit{Characteristic classes
of gauge systems}, Nucl. Phys. \textbf{B703} (2004), 419-453.

\bibitem{M} Mehta R.A., \textit{Supergroupoids, double structures, and equivariant
cohomology}, PhD thesis, University of California, Berkeley, 2006;
\verb"arXiv:math.DG\0605356."


\bibitem{Q} Quillen D., \textit{Superconnections and the Chern character},
Topology \textbf{52}, N1 (1985), 89-95.


\bibitem{Sch}Schwarz~A.S.,
\textit{Semiclassical approximation in Batalin-Vilkovisky
formalizm}, Commun. Math. Phys. \textbf{158} (1993), 373-396.


\bibitem{V} Vaintrob A.Yu., \textit{Lie algebroids and homological vector fields}, Uspekhi Matem. Nauk,  \textbf{52}, N2 (1997), 161-163.

\bibitem{VMP} Voronov A.A., Manin Yu.I., Penkov I.B., \textit{Elements of
supergeometry}, J. Soviet Math. \textbf{51} N1 (1990), 2069-2083.


\end{thebibliography}
\end{document}